\documentclass[12pt]{amsart}
\usepackage{amsmath,amssymb}    
\usepackage{url}                
\usepackage{graphicx}           
\usepackage{amsthm}
\usepackage{amscd}
\usepackage{amsfonts}
\usepackage{latexsym}

\theoremstyle{definition}
\newtheorem{theorem}{Theorem}[section]
\newtheorem{lemma}[theorem]{Lemma}
\newtheorem{proposition}[theorem]{Proposition}
\newtheorem{corollary}[theorem]{Corollary}

\newtheorem{definition}[theorem]{Definition}
\newtheorem{example}[theorem]{Example}

\newtheorem{remark}[theorem]{Remark}

\newcommand{\C}{\mathbb{C}}

\date{}
\newcommand{\h}{{\mathfrak{h}}}

\begin{document}

\pagestyle{plain}
\title{Supports of irreducible spherical representations of 
rational Cherednik algebras of finite Coxeter groups}
\author{Pavel Etingof}
\address{Department of Mathematics, Massachusetts Institute of Technology,
Cambridge, MA 02139, USA}
\email{etingof@math.mit.edu}
\maketitle

{\bf To my father Ilya Etingof on his 80-th birthday, with admiration}

\section{Introduction}

In this paper we determine the support of the irreducible
spherical representation (i.e., the irreducible quotient
of the polynomial representation) of the rational Cherednik
algebra of a finite Coxeter group for any value of the parameter $c$. 
In particular, we determine for which values of $c$ this 
representation is finite dimensional. This generalizes
a result of Varagnolo and Vasserot, \cite{VV}, who classified
finite dimensional spherical representations 
in the case of Weyl groups and equal parameters
(i.e., when $c$ is a constant function). 

Our proof is based on the Macdonald-Mehta integral and the
elementary theory of distributions.  

The organization of the paper is as follows. 
Section 2 contains preliminaries on Coxeter groups 
and Cherednik algebras. In Section 3 we state and prove 
the main result in the case of equal parameters.
In Section 4 we deal with the remaining case of irreducible 
Coxeter groups with two conjugacy classes of reflections. 
In Section 5, as an application, we compute the 
zero set of the kernel of the renormalized Macdonald pairing 
in the trigonometric setting (in the equal parameter case). 
Finally, in the appendix, written by Stephen Griffeth, 
it is shown by a uniform argument (using only the theory of
finite reflection groups) that our classification of 
finite dimensional spherical representations of rational
Cherednik algebras with equal parameters coincides with that of
Varagnolo and Vasserot. 

{\bf Acknowledgements.} It is my great pleasure to dedicate 
this paper to my father Ilya Etingof on his 80-th birthday. 
His selflessness and wisdom made him my main role model, 
and have guided me throughout my life. 
 
This work was  partially supported by the NSF grants
DMS-0504847 and DMS-0854764.

\section{Preliminaries}

\subsection{Coxeter groups}

Let $W$ be a finite Coxeter group of rank $r$ 
with reflection representation $\h_{\Bbb R}$
equipped with a Euclidean $W$-invariant inner product $(,)$.
\footnote{As a basic reference on finite 
Coxeter groups, we use the book \cite{Hu}.} 
Denote by $\h$ the complexification of $\h_{\Bbb R}$. 

For $a\in \h$, let $W_a$ be the stabilizer of $a$ in $W$. 
It is well known that $W_a$ is also a Coxeter group,
with reflection representation $\h/\h^{W_a}\cong (\h^{W_a})^\perp$. 
The group $W_a$ is called a {\it parabolic subgroup} of $W$.  
It is well known that the subgroup 
generated by a subset of the set of simple reflections of $W$
(which corresponds to a subset of nodes of the Dynkin diagram) is
parabolic, and any parabolic subgroup of $W$ is conjugate to one 
of this type. 

Denote by $S$ the set of reflections of $W$. 
For each reflection $s$, pick a vector 
$\alpha_s\in \h_{\Bbb R}$ such that $s\alpha_s=-\alpha_s$  
and $(\alpha_s,\alpha_s)=2$.  
Let 
$$
\Delta_W({\bold x})=\prod_{s\in S}(\alpha_s,{\bold x})
$$ 
be the corresponding {\it discriminant 
polynomial} (it is uniquely determined up to a sign). 

Let $d_i=d_i(W),i=1,...,r$, be the degrees of the
generators of the algebra $\Bbb C[\h]^W$. 
Let $\ell(w)$ be the length of $w\in W$. 
Let 
$$
P_W(q)=\sum_{w\in W}q^{\ell(w)} 
$$
be the Poincar\'e polynomial of $W$. 
Then we have the following well-known identity 
of Bott and Solomon: 
\begin{equation}\label{bs}
P_W(q)=\prod_{i=1}^r \frac{1-q^{d_i}}{1-q}.
\end{equation}

If $W$ is an irreducible Coxeter group which contains 
two conjugacy classes of reflections (i.e. an even dihedral
group $I_2(2m)$ or a Weyl group of type $B_n$, $n\ge 2$, or
$F_4$), then it is useful to consider the 2-variable Poincar\'e
polynomial
$$
P_W(q_1,q_2):=\sum_{w\in W}q_1^{\ell_1(w)}q_2^{\ell_2(w)},
$$
where $\ell_i(w)$ is the number of simple reflections of $i$-th
type occurring in a reduced decomposition of $w$. 

In the Weyl group case, for a positive
root $\alpha$, denote by ${\rm ht}_i(\alpha)$, $i=1,2$, the
number of simple roots of $i$-th type occuring in 
the decomposition of $\alpha$.

\begin{proposition}\label{poincs}(\cite{Ma})
(i) One has 
$$
P_{I_2(2m)}(q_1,q_2)=\frac{1-q_1^2}{1-q_1}\frac{1-q_2^2}{1-q_2}
\frac{1-(q_1q_2)^m}{1-q_1q_2}.
$$

(ii) 
For Weyl groups one has 
$$
P_W(q_1,q_2)=\prod_{\alpha>0}\frac{1-q_\alpha 
q_1^{{\rm ht}_1(\alpha)}
q_2^{{\rm ht}_2(\alpha)}}{1- 
q_1^{{\rm ht}_1(\alpha)}
q_2^{{\rm ht}_2(\alpha)}},
$$
where $q_\alpha=q_i$ if $\alpha$ is a root of $i$-th type, $i=1,2$. 
\end{proposition}

From this proposition one can obtain 
the following more explicit formulas for the 
2-variable Poincar\'e polynomials of $I_2(2m)$, $B_n$ and $F_4$.

\begin{proposition}\label{moreexp} (\cite{Ma}) One has
$$
P_{I_2(2m)}(q_1,q_2)=
(1+q_1)(1+q_2)(1+q_1q_2+...+q_1^{m-1}q_2^{m-1}),
$$
$$
P_{B_n}(q_1,q_2)=\prod_{j=0}^{n-1}(1+q_1+...+q_1^j)
\prod_{j=1}^{n-1}(1+q_1^jq_2),
$$
and
$$
P_{F_4}(q_1,q_2)=
$$
$$
(1+q_1)(1+q_1+q_1^2)(1+q_2)(1+q_2+q_2^2)
(1+q_1^2q_2)(1+q_1q_2^2)(1+q_1q_2)(1+q_1^2q_2^2)(1+q_1^3q_2^3).
$$
\end{proposition}

\subsection{Cherednik algebras}

Let $c$ be a $W$-invariant function on $S$. 
Let $H_c(W,\h)$ be the corresponding rational Cherednik 
algebra (see e.g. \cite{E1}). Namely, $H_c(W,\h)$ 
is the quotient of $\C[W]\ltimes T(\h\oplus \h)$ 
(with the two generating copies of $\h$ spanned by $x_a,y_a$, $a \in \h$), 
by the defining relations
$$
[x_a,x_b]=[y_a,y_b]=0, [y_a,x_b]=(a,b)
-\sum_{s\in S}c_s(\alpha_s,a)(\alpha_s,b)s.
$$

Let $M_c=H_c(W,\h)\otimes_{\Bbb CW\ltimes \Bbb C[y_a]}\Bbb C$, where 
$y_a$ act in $\Bbb C$ by $0$ and $w\in W$ by $1$. 
Then we have a natural vector space isomorphism $M_c\cong \Bbb C[\h]$.
For this reason $M_c$ is called {\em the polynomial
representation} of $H_c(W,\h)$. 
The elements $y_a$ act in this representation by Dunkl operators 
(see \cite{E1}). 

The following proposition is standard, see e.g. \cite{E1}.

\begin{proposition}\label{contra} There exists a unique
$W$-invariant symmetric bilinear form $\beta_c$ on $M_c$
such that $\beta_c(1,1)=1$, which
satisfies the contravariance condition
$$
\beta_c(y_av,v')=\beta_c(v,x_av'),\ v,v'\in M_c, a\in \h.
$$
Polynomials of different degrees are orthogonal under $\beta_c$. 
Moreover, the kernel of $\beta_c$ is the maximal proper submodule 
$J_c$ of $M_c$, so $M_c$ is reducible iff $\beta_c$ is degenerate.  
\end{proposition}

Let $L_c=M_c/J_c$ be the irreducible quotient of $M_c$. 
The module $L_c$ is called {\it the irreducible spherical
representation} of $H_c(W,\h)$. 

\section{The main theorem - the case of equal parameters} 

\subsection{Statement of the theorem}

The goal of this paper is to determine the support of 
$L_c$ as a $\Bbb C[\h]$-module. We start with 
the case when $c$ is a constant function. 
We will assume that $c\in (\Bbb Q\setminus \Bbb Z)_{>0}$; otherwise, 
it is known from \cite{DJO} that $M_c$ is irreducible, so $L_c=M_c$, 
and the support of $L_c$ is the whole space $\h$. 
Let $m$ be the denominator of $c$ (written in lowest terms). 

\begin{theorem}{\label{thm-supp}}
A point $a\in \h$ belongs to the support of $L_{c}$ if and only if
$$
\dfrac{P_{W}}{P_{W_{a}}}(e^{2\pi {\rm i} c})\neq 0, 
$$
i.e., if and only if
$$
\#\{i| m \text{ divides }d_{i}(W)\}=
\#\{i| m \text{ divides }d_{i}(W_a)\}.
$$
\end{theorem}

\begin{remark}
The equivalence of the two conditions in Theorem \ref{thm-supp} follows from 
the Bott-Solomon formula (\ref{bs}) for $P_W$. 
\end{remark} 

\begin{corollary}\label{find}
$L_c$ is finite dimensional if and only if 
$\dfrac{P_{W}}{P_{W'}}(e^{2\pi {\rm i} c})=0$,
i.e., if and only if
$$
\#\{i| m \text{ divides }d_{i}(W)\}>
\#\{i| m \text{ divides }d_{i}(W')\},
$$
for any maximal parabolic subgroup $W'\subset W$.
\end{corollary}

We note that Varagnolo and Vasserot \cite{VV} proved that 
if $W$ is a Weyl group then $L_c$ is finite dimensional if and only 
if there exists a regular elliptic element 
in $W$ of order $m$ (i.e. an element with no eigenvalue $1$ in 
$\h$ and an eigenvector $v$ not fixed
by any reflection, see \cite{Sp}). A direct uniform proof of the 
equivalence of this condition to the condition of Corollary 
\ref{find}, based solely on the theory of finite reflection groups, 
is given in the appendix to this paper, written by S. Griffeth. 

\begin{remark}
If $W$ is a Weyl group, then 
the values of the denominator $m$ of $c>0$ for which $L_c$ is
finite dimensional are listed in \cite{VV}. 
Let us list these values in the noncrystallographic cases. 

For dihedral groups $I_2(p)$: $m\ge 2$ is any number dividing $p$
(this follows from the paper \cite{Chm}). 

For $H_3$: $m=2,6,10$ (this is due to M. Balagovic and
A. Puranik).

For $H_4$: $m$ is any divisor of a degree of $H_4$, i.e. 
$m=2,3,4,5,6,10,12,15,20,30$. 
\end{remark} 

\subsection{Proof of Theorem \ref{thm-supp}} 

\subsubsection{Tempered distributions}

Let $\mathcal  S(\Bbb R^{n})$ be the set of Schwartz functions on 
$\Bbb R^{n}$, i.e. 
$$S(\Bbb R^{n})=\{f\in C^{\infty}(\Bbb R^{n})|
\forall \alpha, \beta, \sup|{\bold x}^{\alpha}\partial^{\beta}f({\bold x})|<\infty\}.
$$
This space has a natural topology.

A tempered distribution on $\Bbb R^{n}$ is a continuous linear functional
on $\mathcal  S(\Bbb R^{n})$. Let $\mathcal  S'(\Bbb R^{n})$ denote the 
space of tempered distributions.

We will use the following well known lemma (see \cite{H}).

\begin{lemma}{\label{lem-dist}}
(i) $\Bbb C[{\bold x}]e^{-{\bold x}^{2}/2}\subset \mathcal  S(\Bbb R^{n})$
is a dense subspace.

(ii) Any tempered distribution $\xi$ has finite order, i.e., $\exists N=N(\xi)$ 
such that if $f\in \mathcal  S(\Bbb R^{n})$ satisfies 
$f={\rm d} f=\cdots={\rm d}^{N-1}f=0$ 
on ${\rm supp} \xi$, then $\langle \xi, f\rangle=0$.
\end{lemma}

\subsubsection{The Macdonald-Mehta integral}

The Macdonald-Mehta integral is the integral 
$$
F_W(c):=(2\pi)^{-r/2}\int_{\h_{\Bbb R}}
e^{-{\bold x}^{2}/2}|\Delta_W({\bold x})|^{-2c}{\rm d} {\bold x}.
$$
It is convergent for ${\rm Re}(c)\le 0$. 

The following theorem gives the value of the Macdonald-Mehta
integral. 

\begin{theorem}\label{mmi}
One has 
$$
F_W(c)=\prod_{i=1}^r \frac{\Gamma(1-d_ic)}{\Gamma(1-c)}.
$$
\end{theorem}

This theorem was conjectured by Macdonald and proved by Opdam \cite{O1}
for Weyl groups and by F. Garvan (using a computer) for $H_3$ and
$H_4$ (for dihedral groups, the formula follows from Euler's beta integral). 
Later, a uniform and computer-free proof for all Coxeter
groups was given in \cite{E2}.

\subsubsection{The Gaussian inner product} 

Let $a_i$ be an orthonormal basis of $\h$, and 
$\bold f=\frac{1}{2}\sum y_{a_i}^2$.
Introduce the {\em Gaussian inner product} on $M_c$ as follows:

\begin{definition}
The Gaussian inner product
$\gamma_c$ on $M_c$
is given by the formula
$$
\gamma_c(v,v')=\beta_c(\exp(\bold f)v,\exp(\bold
f)v').
$$
\end{definition}

This makes sense because the operator $\bold f$ is locally nilpotent on
$M_c$. 

\begin{proposition}\label{intgau} (\cite{Du}, Theorem 3.10)
\footnote{A proof of this theorem can also be found in \cite{ESG}
(Proposition 4.9).}
For $\mathrm{Re}(c) \leq 0 $, one has 
$$
\gamma_{c}(P, Q)=\frac{(2\pi)^{-r/2}}{F_W(c)}\int_{\h_{\Bbb R}}
e^{-{\bold x}^{2}/2}|\Delta_W({\bold x})|^{-2c}P({\bold x})Q({\bold
x}){\rm d} {\bold x}, 
$$
where 
$P, Q$ are polynomials.
\end{proposition}

\subsubsection{Proof of Theorem \ref{thm-supp}}

Consider the distribution
$$
\xi_{c}^{W}=\frac{(2\pi)^{-r/2}}{F_W(c)}|\Delta_W({\bold x})|^{-2c}.
$$
It is well-known that this distribution extends to a meromorphic 
distribution in $c$ (Bernstein's theorem). Moreover, since 
$\gamma_c(P,Q)$ is a polynomial in $c$ for any $P$ and $Q$,  
this distribution is in fact holomorphic in $c\in \Bbb C$.

\begin{proposition}{\label{prop-1}}
\begin{eqnarray*}
{\rm supp}(\xi_{c}^{W})&=&
\{a\in \h_{\Bbb R}|\frac{F_{W_{a}}}{F_{W}}(c)\neq 0\}
=\{a\in \h_{\Bbb R}|\frac{P_{W}}{P_{W_{a}}}(e^{2\pi\mathrm{i}c})\neq 0\}\\
&=&\{a\in \h_{\Bbb R}|\#\{i|\text{denominator of $c$ divides } d_{i}(W)\}\\
&&\qquad\qquad=\#\{i|\text{denominator of $c$ divides } d_{i}(W_{a})\}\}.
\end{eqnarray*}
\end{proposition}

\begin{proof}
First note that the last equality 
follows from the Bott-Solomon formula (\ref{bs}) 
for the Poincar\'e polynomial, and 
the second equality from 
Theorem \ref{mmi}. Now let us prove the first equality. 

Look at $\xi_{c}^{W}$ near $a\in \h$. Equivalently, 
we can consider 
$$\xi^{W}_{c}({\bold x}+a)=\frac{(2\pi)^{-r/2}}{F_W(c)}|\Delta({\bold x}+a)|^{-2c}$$
with ${\bold x}$ near $0$.
We have
\begin{eqnarray*}
\Delta_{W}({\bold x}+a)&=&\prod_{s\in S}\alpha_{s}({\bold x}+a)
=\prod_{s\in S}(\alpha_{s}({\bold x})+\alpha_{s}(a))\\
&=&\prod_{s\in S\cap W_{a}}\alpha_{s}({\bold x})\cdot
\prod_{s\in S\backslash S\cap W_{a}}(\alpha_{s}({\bold x})+\alpha_{s}(a))\\
&=&\Delta_{W_{a}}({\bold x})\cdot G({\bold x}),
\end{eqnarray*}
where $G$ is a nonvanishing function near $0$
(since $\alpha_s(a)\ne 0$ if $s\notin S\cap W_a$).

So near $0$, we have
$$
\xi_{c}^{W}({\bold x}+a)=
\frac{F_{W_{a}}}{F_{W}}(c)\cdot \xi_{c}^{W_{a}}({\bold x})\cdot 
|G(\bold x)|^{-2c},
$$
and the last factor is well defined since $G$ is nonvanishing.
Thus $\xi_{c}^{W}({\bold x})$ is nonzero near $a$ if and only if
$\frac{F_{W_{a}}}{F_{W}}(c)\neq 0$,
which finishes the proof.
\end{proof}

Now consider the support of $L_c$. 
Note that $\h$ has a stratification by stabilizers of points in $W$, and
by the results of \cite{Gi} (see also \cite{BE}), 
the support of $L_c$ is a union of strata of this
stratification. 

\begin{proposition}{\label{prop-2}} For any $c\in \Bbb C$,
$${\rm supp}(\xi_{c}^{W})=({\rm supp} L_{c})_{\Bbb R},$$
where the right hand side denotes the set of 
real points of the support. 
\end{proposition}

\begin{proof}
Let $a\notin {\rm supp} L_{c}$ 
and assume $a\in {\rm supp} \xi^{W}_{c}$.
Then we can find a $P\in J_{c}=\ker \gamma_{c}$ 
such that 
$P(a)\neq 0$.
Pick a compactly supported test function $\phi\in C_{c}^{\infty}(\h_{\Bbb{R}})$ such that 
$P$ does not vanish anywhere on ${\rm supp} \phi$, and $\langle\xi^{W}_{c},\phi\rangle\neq 0$
(this can be done since $P(a)\ne 0$ and $\xi^W_c$ is nonzero near $a$). Then we have
$\phi/P\in \mathcal {S}(\h_{\Bbb R})$.
Thus from Lemma \ref{lem-dist}(i) it follows that 
there exists a sequence of polynomials $P_{n}$ such that
$$P_{n}({\bold x})e^{-{\bold x}^{2}/2}\to \frac{\phi}{P} \text{ in }
\mathcal {S}(\h_{\Bbb R})\text{, when }n\to \infty.$$
So $PP_{n}e^{-{\bold x}^{2}/2}\to \phi\text{ in }
\mathcal {S}(\h_{\Bbb R})\text{, when }n\to \infty.$

But by Proposition \ref{intgau}, we have 
$\langle \xi^{W}_{c},PP_{n}e^{-{\bold x}^{2}/2}\rangle=
\gamma_{c}(P, P_{n})$. Hence, 
$\langle \xi^{W}_{c},PP_{n}e^{-{\bold x}^{2}/2}\rangle=0$,
which is a contradiction. This implies that 
${\rm supp} \xi_c^W\subset ({\rm supp} L_c)_{\Bbb R}$. 

To establish the opposite inclusion, 
let $P$ be a polynomial on $\h$ which vanishes identically on 
${\rm supp} \xi_c^W$. By Lemma \ref{lem-dist}(ii), there exists $N$ such that 
$\langle \xi_c^W, P^N({\bold x})Q({\bold x})e^{-{\bold x}^2/2}\rangle =0$. 
Thus, using Proposition \ref{intgau}, we see that for any 
polynomial $Q$, $\gamma_c(P^N,Q)=0$, i.e. $P^N\in {\rm Ker}\gamma_c$. 
Thus, $P|_{{\rm supp} L_c}=0$. This implies the required inclusion, since
${\rm supp}  \xi_c^W$ is a union of strata.  
\end{proof}

Theorem \ref{thm-supp} follows 
from Proposition \ref{prop-1} and Proposition \ref{prop-2}.

\section{The main theorem - the case of non-equal parameters} 

\subsection{Statement of the theorem} 

Consider now the case when $W$ is an irreducible Coxeter group
with two conjugacy classes of reflections.
In this case, $c=(c_1,c_2)$,
and by $e^{2\pi i c}$ we will mean the pair
$(q_1,q_2)$, where $q_j=e^{2\pi ic_j}$, $j=1,2$. 

Define a {\it positive} line in the plane with coordinates 
$(c_1,c_2)$ to be any line of the form $a_1c_1+a_2c_2=b$, where 
$a_1,a_2\ge 0, b > 0$. 

\begin{theorem}\label{thsupp2}
A point $a\in \h$ belongs to the support 
of $L_c$ if and only if there is no 
positive line passing through $c$ on which the function 
$z\mapsto \frac{P_W}{P_{W_a}}(e^{2\pi iz})$ identically 
vanishes. 
\end{theorem}

\begin{corollary}\label{findi}
$L_c$ is finite dimensional if and only if 
for every maximal parabolic subgroup $W'\subset W$, 
there exists a positive line $\ell$ passing through $c$ such that
the function $\frac{P_W}{P_{W'}}(e^{2\pi iz})$ vanishes on $\ell$.  
\end{corollary} 
 
\subsection{Computation of points $c$ for which $L_c$ is finite dimensional}

Let us use Corollary \ref{findi} to compute explicitly the set
$\Sigma_c$ of points $c$ for which $L_c$ is finite dimensional. 
The computation is straightforward using Propositions 
\ref{poincs} and \ref{moreexp} (although somewhat tedious),
so we will only give the result. 

{\bf 1. The dihedral case, $I_2(2m)$.} In this case, the set
$\Sigma_c$ is the union of the following lines and isolated points.

1) The lines are $c_1+c_2=\frac{r}{m}$, where $r\in \Bbb N$ and $r$ is
not divisible by $m$.

2) The isolated points are $(\frac{p_1}{2},\frac{p_2}{2})$,  
where $p_j$ are odd positive integers. 

This description coincides with the one of \cite{Chm}. 
 
{\bf 2. The case $F_4$.}
In this case, the set
$\Sigma_c$ is the union of
the following lines and isolated points.

1) The lines are $c_1+c_2=\frac{p}{4}$ and
$c_1+c_2=\frac{p}{6}$, where $p$ is an odd positive integer. 

2) The isolated points are: 

2a: $(\frac{p_1}{2},\frac{p_2}{2})$, where $p_1,p_2$ are
odd positive integers; 

2b: $(\frac{p_1}{3},\frac{p_2}{3})$, 
where $p_1,p_2$ are positive integers not divisible by $3$; 

2c: $(\frac{p_1}{3},\frac{p_2}{4}-\frac{p_1}{6})$
and $(\frac{p_2}{4}-\frac{p_1}{6},\frac{p_1}{3})$, 
where $p_2$ is an odd positive integer and $p_1$ is a positive
integer not divisible by $3$;

2d: $(\frac{2p_2-p_1}{6},\frac{2p_1-p_2}{6})$, where $p_1,p_2$ 
are odd positive integers such that $p_1+p_2$ is not divisible by $3$.
 
{\bf 3. The case $B_n$, $n\ge 2$.}
In this case, let $c_1$ correspond to long roots. Then the set
$\Sigma_c$ is the union of
the following lines and isolated points.

1) The lines are $(n-1)c_1+c_2=\frac{p}{2}$,
where $p$ is an odd positive integer. 

2) The isolated points are
$(\frac{r}{n},\frac{p}{2}-r+\frac{rs}{n})$, 
where $r$ is a positive integer not divisible by $n$, 
$p$ is an odd positive integer, and $2\le s\le \frac{n}{{\rm gcd}(r,n)}$ is an integer.

\begin{remark}
Note that in the case $c_1=c_2$, 
we recover precisely the result of Varagnolo and Vasserot,
\cite{VV} for $W$ of types $B_n$, $F_4$ and $G_2$,
while setting $c_2=0$, we recover their result for $W$ of type $D_n$.  
\end{remark}

\subsection{Proof of Theorem \ref{thsupp2}}

First we need to formulate the 
appropriate generalization of the
Macdonald-Mehta integral. 
Let $S_1,S_2$ be the sets of reflections in $W$ of the first and
second kind, and let 
$$
\Delta_{W,j}(\bold x)=\prod_{s\in S_j}(\alpha_s,\bold x).
$$
Define the Macdonald-Mehta integral with two parameters: 

$$
F_W(c_1,c_2):=(2\pi)^{-r/2}\int_{\h_{\Bbb R}}
e^{-{\bold x}^{2}/2}|\Delta_{W,1}({\bold x})|^{-2c_1}
|\Delta_{W,2}({\bold x})|^{-2c_2}{\rm d} {\bold x}.
$$
As before, it is convergent for ${\rm Re}(c_j)\le 0$. 

The following theorem gives the value of the two-parameter Macdonald-Mehta
integral. 

\begin{theorem}\label{mmi2} 
(i) For dihedral groups $I_2(2m)$, one has 
$$
F_W(c_1,c_2)=
\frac{\Gamma(1-2c_1)}{\Gamma(1-c_1)}
\frac{\Gamma(1-2c_2)}{\Gamma(1-c_2)}
\frac{\Gamma(1-m(c_1+c_2))}{\Gamma(1-(c_1+c_2))}.
$$

(ii) For Weyl groups, one has 
$$
F_W(c_1,c_2)=\prod_{\alpha>0}
\frac{\Gamma(1-c_\alpha-c_1{\rm
ht}_1(\alpha)-c_2{\rm ht}_2(\alpha))}
{\Gamma(1-c_1{\rm
ht}_1(\alpha)-c_2{\rm ht}_2(\alpha))}.
$$
\end{theorem}

\begin{proof} (i) follows from Euler's beta integral, and 
(ii) is proved in \cite{O1}. 
\end{proof}

Also, we need an analog of the integral formula for the
Gaussian form $\gamma_c$. This analog  
is given by the following proposition,
whose proof is a straightforward generalization of 
the proof of Proposition \ref{intgau}.  

\begin{proposition}\label{intgau1} 
For $\mathrm{Re}(c_j) \leq 0 $, one has 
$$
\gamma_{c}(P, Q)=\frac{(2\pi)^{-r/2}}{F_W(c_1,c_2)}\int_{\h_{\Bbb R}}
e^{-{\bold x}^{2}/2}|\Delta_{W,1}({\bold x})|^{-2c_1}
|\Delta_{W,2}({\bold x})|^{-2c_2}P({\bold x})Q({\bold
x}){\rm d} {\bold x}, 
$$
where 
$P, Q$ are polynomials.
\end{proposition}

Now we are ready to prove Theorem \ref{thsupp2}. 
Consider the distribution
$$
\xi_{c}^{W}=\frac{(2\pi)^{-r/2}}{F_W(c_1,c_2)}|\Delta_{W,1}({\bold
x})|^{-2c_1}|\Delta_{W,2}({\bold x})|^{-2c_2}.
$$
As before, this distribution extends to a meromorphic 
distribution in $c$ (by Bernstein's theorem), and since 
$\gamma_c(P,Q)$ is a polynomial in $c$ for any $P$ and $Q$,  
this distribution is in fact holomorphic in $c$.

\begin{proposition}{\label{prop-11}}
One has
$$
{\rm supp}(\xi_{c}^{W})=
\{a\in \h_{\Bbb R}|\frac{F_{W_{a}}}{F_{W}}(c)\neq 0\}.
$$
\end{proposition}

\begin{proof}
The proof is parallel to the proof of 
Proposition \ref{prop-1}.
\end{proof}

\begin{corollary}\label{coro}
A point $a\in \h_{\Bbb R}$ belongs to the support 
of $\xi_c^W$ if and only if there is no 
positive line passing through $c$ on which the function 
$z\mapsto \frac{P_W}{P_{W_a}}(e^{2\pi iz})$ identically vanishes. 
\end{corollary} 

\begin{proof}
The Corollary follows from Propositions \ref{prop-11} and
\ref{poincs} and Theorem
\ref{mmi2}, using the bijective correspondence between 
the factors in the product formulas
for $P_W(q_1,q_2)$ in Proposition \ref{poincs} and the
$\Gamma$-factors in the product formulas for $F_W(c_1,c_2)$
in Theorem \ref{mmi2}.   
\end{proof}

\begin{proposition}{\label{prop-21}} For any $c\in \Bbb C^2$,
$${\rm supp}(\xi_{c}^{W})=({\rm supp} L_{c})_{\Bbb R}.$$ 
\end{proposition}

\begin{proof}
Parallel to the proof of Proposition \ref{prop-2},
using Proposition \ref{intgau1}.
\end{proof}

Proposition \ref{prop-11} and \ref{prop-21} imply Theorem \ref{thsupp2}.

\section{Application: the zero set of the kernel of the renormalized Macdonald pairing} 

As an application of the above technique, let us compute the 
zero set of the kernel of the renormalized Macdonald pairing 
in the trigonometric setting (in the equal parameter case). 

Let $R$ be a reduced irreducible root system of rank 
$r$, $W$ be the Weyl group of $R$, $R_+$ a system of positive
roots, $P$ the weight lattice of $R$, and $H={\rm Hom}(P,\Bbb C^*)$ 
be the complex torus attached to $R$. 
Let $\h={\rm Lie}(H)$. For $a\in H$, the stabilizer $W_a$ 
is a reflection subgroup of $W$ which is not necessarily parabolic; 
such a subgroup is called a quasiparabolic subgroup. 

The following lemma is well known. 

\begin{lemma} If $W'\subset W$ is any reflection subgroup, then 
$P_W(q)$ is divisible by $P_{W'}(q)$. 
\end{lemma}

\begin{proof} $P_W/P_{W'}$ is the character of the generators 
of the module $\Bbb C[\h]^{W'}$ over $\Bbb C[\h]^W$. 
\end{proof}

Let 
$$
D_R=\prod_{\alpha\in R_+}(e^\alpha-1)
$$
be the Weyl denominator of $R$. Let $\Bbb C[H]$ 
denote the algebra of regular functions on $H$. 
Let $c\in \Bbb C$. Recall that the {\it Macdonald pairing} 
on the space $\Bbb C[H]$ (or $\Bbb C[H]^W$) is defined by the formula 
$$
\langle P,Q\rangle_c:=\int_{H_{\Bbb R}} |D_R(\bold t)|^{-2c}P(\bold t)Q(\bold t)d\bold t,
$$
where $H_{\Bbb R}$ is the maximal compact subgroup of $H$. 
It is known that this pairing is well defined when ${\rm Re}(c)< 0$, 
and develops poles when $c>0$ and $e^{2\pi ic}$ is a root of $P_W$. 

For $c\in \Bbb C$, ler $d$ be the order of the pole of the 
Macdonald pairing at $c$. Define the {\it renormalized Macdonald pairing} 
by the formula  
$$
(P,Q)_{c}=\lim_{k\to c}(k-c)^d\langle P,Q\rangle_k.
$$
This pairing is well defined and nonzero for any $c$.
 
Moreover, it is easy to see that the kernel $I_c$
of this pairing is an ideal in $\Bbb C[H]$; in fact, 
this ideal is invariant under the trigonometric Dunkl operators, 
so the quotient ring $V_c:=\Bbb C[H]/I_c$ is a representation of 
the trigonometric Cherednik algebra $H_c^{\rm trig}(W,H)$ (i.e., degenerate double
affine Hecke algebra, see \cite{Ch}, \cite{EM}). Our goal 
is to find the support ${\rm supp}(V_c)$ as a $\Bbb C[H]$-module, 
i.e. the zero set of $I_c$. 

The main result of this section is the following theorem.

\begin{theorem}\label{trigg} Suppose that $c\in \Bbb Q_{>0}$. 
A point $a\in H$ belongs to ${\rm supp}(V_c)$ if and only if
the polynomial $P_W/P_{W_a}$ takes a nonzero value at $e^{2\pi i c}$. 
\end{theorem}

\begin{proof}
We have a stratification of $H$ by stabilizers of points in $W$. 
By the results of \cite{BE}, ${\rm supp}(V_c)$ is a union of strata of 
this stratification. Thus, it suffices to prove the result
for $a\in H_{\Bbb R}$.

Consider the distribution on $H_{\Bbb R}$ defined by the formula
$$
\xi_c^R:=\frac{(2\pi)^{-r/2}}{F_W(c)}|D_R|^{-2c}.
$$
This distribution is defined for all $c$, and 
up to a scalar, 
$$
(P,Q)_c=\langle \xi_c^R,PQ\rangle.
$$

Let ${\rm supp}(V_c)_{\Bbb R}$ 
be the intersection of ${\rm supp}(V_c)$ with $H_{\Bbb R}$. 

\begin{proposition}\label{l1}
A point $a\in H_{\Bbb R}$ belongs to ${\rm supp}(\xi_c^R)$ if and only if
the polynomial $P_W/P_{W_a}$ takes a nonzero value at $e^{2\pi i c}$. 
\end{proposition} 

\begin{proof}
The proof is parallel to the proof of Proposition \ref{prop-1}.
Namely, for small $\bold x\in \h_{\Bbb R}$, we have 
$$
\xi_c^R(ae^{\bold x})=\frac{F_{W_a}(c)}{F_W(c)}
\xi_c^{W_a}(\bold x)f(\bold x),
$$ 
where $f$ is a nonvanishing smooth function. 
So the result follows. 
\end{proof}

\begin{proposition}\label{l2}
The set ${\rm supp}(V_c)_{\Bbb R}$ coincides with the support 
${\rm supp}(\xi_c^R)$ of the distribution $\xi_c^R$.  
\end{proposition}

\begin{proof}
The proof is parallel to the proof of Proposition \ref{prop-2}.
\end{proof} 

Theorem \ref{trigg} follows from Proposition \ref{prop-1} and
Proposition \ref{prop-2}.
\end{proof}

\begin{example}
Let $R$ be the root system of type $B_2$. Then 
$P_W(q)=(1+q)^2(1+q^2)$, so setting $q=e^{2\pi ic}$, we see 
that at $c=1/2$ (i.e. $q=-1$), the only points 
$a\in H=(\Bbb C^*)^2$ for which $P_W/P_{W_a}$ does not vanish at $q$
are $(1,1)$ ($W_a=W$) and $(-1,-1)$ ($W_a=\Bbb Z_2\times
\Bbb Z_2$, $P_{W_a}=(1+q)^2$). So the support of $V_c$ is the set
consisting of these two points. 

This example shows that unlike the rational case, the module
$V_c$ is not necessarily irreducible. Namely, local analysis 
near the two points of the support using the results of
\cite{Chm} shows that $V_c$ is the direct
sum of a 1-dimensional irreducible representation supported at
$(-1,-1)$ and a 4-dimensional irreducible representation
supported at $(1,1)$. 
\end{example}

\section{Appendix}

\vskip .05in
\centerline{\bf by Stephen Griffeth}
\vskip .05in

Let $W$ be
a finite real reflection group with reflection representation $\h$.
Recall that an \emph{elliptic element} of $W$ is an element not
contained in any proper parabolic subgroup, 
or, equivalently, an element with fix space $\{0 \}$ in $\h$. 
Recall also that a positive integer $m$ is a
\emph{regular number} for $W$ if there is an element $g \in W$ 
that has a regular eigenvector (i.e., one not fixed by any
reflection) with eigenvalue a primitive $m$-th root of $1$.  
(By Theorem 4.2(i) of \cite{Sp}, in this case the order
of $g$ is $m$; such elements are called \emph{regular}). 
If in addition this element $g$ 
can be chosen to be elliptic, then $m$ is called an 
\emph{elliptic number} for $W$. 

Let $d_1(W),\dots,d_r(W)$ be the degrees of $W$, and let $m$ be a
positive integer. Denote by $a_W(m)$ the number 
of degrees divisible by $m$: 
$a_W(m)=\# \{1 \leq i \leq r \ | \ m \ \text{divides} \ d_i(W)
\}$.  

The purpose of this appendix is to give a uniform proof of the
following theorem. 

\begin{theorem}
Let $W$ be a finite real reflection group. Then $m$ is an
elliptic number for $W$ if and only if for
every maximal parabolic subgroup $W'$ of $W$, 
one has $a_W(m)>a_{W'}(m)$. 
\end{theorem}

\begin{proof}
First suppose $m$ is an elliptic number for $W$.
This means that there exists an elliptic element $b\in W$ 
and a regular vector $v\in \h$ 
such that $bv=\zeta v$, where $\zeta$ 
is a primitive $m$-th root of unity. Assume towards a
contradiction that $a_W(m)\le a_{W'}(m)$ for some maximal
parabolic subgroup $W'$. Then by part (i) of Theorem 3.4 of \cite{Sp},
$a_W(m)=a_{W'}(m)$, and there is an element $g \in W'$ so that the 
$\zeta$-eigenspace of $g$ has dimension
exactly equal to $a_W(m)$. By part (iv) of Theorem 4.2 of
\cite{Sp}, the elements $b$ and $g$ are conjugate in $W$. 
This is a contradiction, since $b$ is an elliptic elememt, and
$g$ is not.  

Conversely, assume the inequalities in the statement of the theorem
hold.  These inequalities together with Part (i) of Theorem 3.4 of
\cite{Sp} imply that for any primitive $m$-th root of unity $\zeta$ 
there exists an element $g \in W$ with $\zeta$-eigenspace of 
dimension $a_W(m)$,
and the fix space of any such $g$ in $\h$ is zero 
(i.e., $g$ is elliptic). Since $W$ is a real reflection
group, this implies that the determinant of $g$ on $\h$ 
is $(-1)^r$ (i.e., is independent of $g$).  Examining
the left hand side of the equation in Corollary 2.6 of \cite{LM} shows
that the term $(-T)^{a_W(m)}$ occurs with non-zero coefficient. 
Hence, looking at the right hand side of this equation, 
we see that the number of codegrees of $W$ divisible by $m$ is
$a_W(m)$. Now part (ii) of Theorem 3.1 of \cite{LM} 
implies that $m$ is a regular number, and
hence elliptic.
\end{proof}

{\bf Acknowledgements.} 
Stephen Griffeth acknowledges the full financial support of EPSRC grant EP/G007632.

\end{document}